# C.S. Peirce's Classification of Dyadic Relations:
## Exploring the Relevance for Logic and Mathematics

Jeffrey Downard[1]

In "The Logic of Mathematics; an Attempt to Develop my Categories from Within" [CP 1.417-520] Charles S. Peirce develops a scheme for classifying different kinds of monadic, dyadic and triadic relations. His account of these different classes of relations figures prominently in many parts of his philosophical system, including the phenomenological account of the categories of experience, the semiotic account of the relations between signs, objects and the metaphysical explanations of the nature of such things as chance, brute existence, law-governed regularities and the making and breaking of habits. Our aim in this essay is to reconstruct and examine central features of the classificatory system that he develops in this essay. Given the complexity of the system, we will focus our attention in this essay mainly on different classes of degenerate and genuine dyadic relations, and we will take up the discussion of triadic relations in a companion piece.

One of our reasons for wanting to explore Peirce's philosophical account of relations is to better understand how it might have informed the later development of relations as they figure in the history of mathematical logic. Peirce's writings on formal relations had a direct and significant impact on the Ernest Schröder's development of his dyadic logical systems. Together, the earlier works of Peirce and Schröder influenced the development of the seminal works of Leopold Löwenheim, Thoralf Skolem and Alfred Tarski [Anellis 2004].

---

[1] Associate Professor, Department of Philosophy, Northern Arizona University.

[2] In what follows, I will refer to Peirce's essay "The Logic of Mathematics; An Attempt to Develop my Categories from Within" using the abbreviated title "The Logic of Mathematics."
[3] One reason to retain all of the parts of the figures—including the diamonds that represent the saturated bonds—is that it makes it possible to study more closely the combinatorial possibilities of the system.

[4] The operation of union in virtue of which such a set is formed can help us understand the relationship between the number two considered cardinally and the number two considered ordinarily.
[5] This account of unity and identity does not appear to be subject to any of the objections that Frege raises against alternate accounts drawn from the history of mathematics and philosophy. Furthermore, Peirce's explanation of the character of these conceptions seems considerably richer

The essay has six parts. First, we will provide an overview of the main lines of argument in "The Logic of Mathematics".[2] Second, we will examine his account of the material categories and consider the relationship of these categories to the formal elements in experience. Third, we will clarify the points Peirce is making about the formal character of monadic, dyadic and triadic relations by developing a graphical system for diagramming the different kinds of structures that can be composed by combining these elementary relations. Fourth, we will apply this graphical system to the analyses of the main classes of formal relations in "The Logic of Mathematics" with special emphasis on the classification of dyadic relations. For the purposes of this inquiry, we will draw on three main sources: math, phenomenology and logic.

## 1. An overview of "The Logic of Mathematics"

In his seminal exposition of the main themes in Peirce's work in philosophy, logic and mathematics, Murray Murphey raises the following concerns about the strategy Peirce has adopted for answering the three questions about mathematical hypotheses [Murphey 1961]. His concerns are directed towards Peirce's account of the relationship between the mathematics, on the one hand, and the logical and phenomenological account of the categories, on the other hand. First, on the basis of architectonic considerations, he suggests that the purported connection between the categories and mathematics tends to introduce "confusion into the theory of the categories." What is more, Murphey says:

> Since Peirce's unwavering allegiance to the architectonic plan required that numbers too should fall within the categorical schema, they are classed as Firsts. But since Firsts are represented by icons, and, as we have seen above, icons also serve as variables, Peirce's dictum that mathematical propositions are iconic now becomes ambiguous in an almost vicious sense…. Here as in so many other cases, the categorical schema fails to preserve significant distinctions and so helps to obscure what it is that Peirce is talking about [Murphey 1961, pp. 239-40].

---

[2] In what follows, I will refer to Peirce's essay "The Logic of Mathematics; An Attempt to Develop my Categories from Within" using the abbreviated title "The Logic of Mathematics."

In this essay, we seek to show that these kinds of concerns are misplaced. On the one hand, Peirce is quite clear in claiming that the signs we use to refer to numbers in processes of counting are, first and foremost, indexes and not icons or symbols. Having said that, we should not lose sight of the fact that every index is able to refer to its object only in virtue of the icons that pick out the properties of the objects [CP 4.544]. In the case of numbers, the icons that seem to be of primary importance are not the qualisigns that pick out the color of the ink used to scribe the figures or the color of the paper upon which those figures are written. Certainly, our ability to observe the mathematical relations does depend in some sense upon our ability to discriminate between the color of the ink and the color of the paper. But the observation of the mathematical relations involves an abstractive process that is precisive in character in which we attend to the formal relations and ignore the material characters that are largely irrelevant to the processes involved in cognizing the mathematical relations [CP 4.234-5].

Murphey raises a second set of concerns that are based on worries about the coherence of Peirce's phenomenological account of the categories. He states his worries in the following way:

> It is impossible to regard Peirce's phenomenological treatment of the categories as anything more than a quite unsuccessful sleight of hand. Even if his attempt to identify formal logic with mathematics is accepted, the most that results from it is that there is an algebra isomorphic to the logical system of relations. Now this algebra could certainly be used to classify elements of the phaneron, if those elements should happen to exhibit characteristics which would admit of such a classification, and it is the purpose of the phenomenology to show that in fact such a classification can be made. But what the phenomenology does not show is why it should be made. There are certainly other ways of classifying the elements of the phaneron which are equally simply and exhaustive, and no reason is given as to why the classification by relations is to be preferred. [Murphey 1961, pg. 368].

In this essay, we will attempt to determine whether or not Murphey's concerns about the coherence of Peirce's phenomenological account of the categories are on track or misguided. I will argue that, in fact, the concerns are misplaced because they rest on mistaken assumptions about the purposes for making a phenomenological classification of the categorical elements in ordinary experience and about the relationship between what can be learned from a logical study of reasoning in mathematical inquiry and the phenomenological analysis of the formal elements.

Before we dig into the details of the essay, I think is appropriate to start with an overview of the larger argument. My hope is that this quick overview will provide some reasons for believing that we should take a closer look at the questions Peirce is trying to answer and the way he is drawing on the phenomenological and the logical accounts of the categories in order to develop and refine his answers.

One reason the essay is puzzling and can easily lead a reader to misinterpret what Peirce is saying about the relationships that hold between mathematics, phenomenology and logic is that several of the steps in the argument are suppressed—especially those that form the initial premises. One reason the initial steps are suppressed is that the first four pages of the manuscripts from which the published version is drawn are missing from the archival collection [CP 1.417 n1]. As a result, we will need to make some educated guesses about how Peirce is framing the questions we are taking as the starting points for our inquiries.

For starters, let us make some of the implicit steps in the argument more explicit by drawing them out of related essays where Peirce does make the points more clearly. Let's draw out the implicit premises in a manner that will help to clarify the methods of inquiry that are being used to study (1) the logic of mathematics, (2) the phenomenological account of the categories and (3) the logical and semiotic account of the categories. In this essay from 1896, Peirce had not yet made the kind of clear separation between the phenomenological and logical accounts of the categories that he later makes starting around 1903 when he decides to treat phenomenology and semiotics as a separate branches of inquiry with their own aims, methods and types of observations. [Peirce 1998, pp. 145-59 and 360-70]. My intention in presenting a reconstruction of the argument of "The Logic of Mathematics" is to read it in light of these later developments in Peirce's understanding of how the phenomenological and logical inquiries might be better fitted together.

Let us try to summarize the major moves in the argument. First, every area of inquiry, including those in the special sciences as well as those in the formal sciences of mathematics and the cenoscopic sciences of philosophy, must start with observations. What is more, every area of inquiry should draw on some version of the experimental method—specially adapted to the kinds of observations that are at hand—in order to develop more systematic explanations of the phenomena that call out for explanation. The experimental method is a form of inquiry that typically involves the following sorts of steps applied in an iterative fashion to correct for errors and improve the security and the uberty of the conclusions: (1) start with accepted theories; (2) observe a surprising phenomenon (3) formulate competing hypotheses by abduction to explain the phenomenon; (4) draw out the consequences by deduction from the competing hypotheses; (5) conduct experiments to test the rival hypotheses and determine by induction what is and is not supported by the data [W 3.323-38].

Second, philosophical inquiry must start from observations that are drawn from common experience [CP 3.428]. This "common experience" includes all that *could* be felt, thought or experienced in any way by any person who is awake at any moment of any day. The various features of this common experience provide us with the observations that serve as the starting points of philosophical inquiry. The simple fact that this experience is so familiar to us is the reason that philosophers seem to have particular difficulty in seeing what is staring them in the face [CP 1.134].

Third, the purpose of phenomenological inquiry is to put us in a better position to analyze in more minute detail what is contained in the observations that are being made of the various sorts of phenomena that seem to call out of explanation, and then to correct for any observational errors that have occurred [Downard 2014]. Instead of waiting until hypotheses are tested, we can make considerable progress in inquiry by more closely examining the phenomena that have been observed in order to identify possible sources of observational error.

Fourth, in order to make more penetrating analyses of the various sorts of observations that can be drawn from common experience, we need an account of the elements found in these experiences. On Peirce's account, there are three basic material categories that are found in all experience: quality, brute fact and thought. What is more, there are three formal elements that are necessary features of all possible experience, and these formal elements correspond with the three material categories in the following way: the material category of quality corresponds with the formal category of the *monad* or what has firstness; the material category of brute fact corresponds with the formal category of the *dyad* or what has secondness; the material category of thought or lawfulness corresponds with the formal category of the *triad* or what has thirdness [CP 1.417-21].

Fifth, the numerical conceptions of one, two, and three are intimately related to the monad, dyad, and triad. The formal features found in all experience are elements that we abstract in a precise manner from the observations in order to formulate the conceptions that comprise the most universal of the hypotheses that lie at the bases of mathematics [CP 1.441-2 and 1.471].

The upshot of the argument is that mathematics, like philosophy, is an observational science. Unlike the empirical sciences, mathematics and philosophy largely draw on the observations of common experience. They do not require special observations made with telescopes or microscopes. Unlike philosophy or the special sciences, however, mathematics is not a positive science. Rather, it is an entirely formal science. That is, mathematics does not attempt to establish positive matters of fact about what really is the case concerning the nature of things. Instead, it studies the formal relations between the parts of ideal systems [CP 3.558]. In the context of providing a detailed classification of the sciences, he says:

> The first is mathematics, which does not undertake to ascertain any matter of fact whatever, but merely posits hypotheses, and traces out their consequences. It is observational, in so far as it makes constructions in the imagination according to abstract precepts, and then observes these imaginary objects, finding in them relations of parts not specified in the precept of construction. This is truly observation, yet certainly in a very peculiar sense; and no other kind of observation would at all answer the purpose of mathematics [CP 1.240].

In what follows, our aim is to fill in more of the details of this general line of argument by drawing out the role of the categories in constructing mathematical hypotheses, formulating diagrams and then reasoning deductively from hypotheses and diagrams.

## 2. A phenomenological account of the categories of experience

Peirce argues that the formal relations of the monad, dyad, and triad are found in every part of our ordinary experience, and they have the same basic structure as the relations that are essential for constructing idealized figures and reasoning about mathematical diagrams. This should not come as too much of a surprise given the fact that mathematical diagrams are abstractions formed on the basis of everyday experience.

Taken in isolation, each of the categories is but one element in our experience. The first of the three material categories of experience is comprised of the qualities that we feel in any phenomena we might observe. The examples of qualities that he considers include the qualities of red, bitter, tedious, hard, heartrending and noble [CP 1.422]. In each case, we could say of some thing, such as a fire engine, that it "is red," or of a cup of coffee that it "is bitter." In doing so, we are abstracting the quality from the object of which it is a characteristic. The type of abstraction is precisive in character, such that we attend to the redness or the bitterness and ignore the fire engine or the cup of coffee [CP 4.463]. In doing this, we make a logical separation of a particular type between the representation of the object as a thing that is located at some time and place and the feeling of the quality insofar as it is considered to be sign that functions in our judgments and assertions. The logical part of the analysis enables us to see that the conception of the qualitative characteristic has the character of a monadic predicate. The phenomenological part of the analysis enables us to see that the quality, considered as a possible feeling, is something that can be similar to other qualities, such that red is similar to pink and bitter is similar to sour. While we often make comparisons in our ordinary experience between qualities in terms of both their similarities and differences, the aim of the phenomenological analysis is to separate the material element from its relations to any parts that we might take it to have, such as the intensity of the experience, or to any relations that a given quality might have to other qualities. A quality, considered in complete abstraction from all such considerations, has the character of what is "eternal, independent of time and of any realization" [CP 1.420]. As such, the pure qualities have no individual identities. In this respect, the experience of a feeling of quality is strikingly different from the experience of an object that is here and now.

    Empiricist philosophers have raised metaphysical questions about the nature of qualities, with Locke arguing that the primary properties of objects are different in kind from the secondary qualities of our feelings [Locke 1975]. For the methodological purposes governing phenomenological inquiry, however, Peirce argues that we need to

set the metaphysical questions to the side. It is sufficient that "wherever there is a phenomenon there is a quality; so that it might almost seem that there is nothing else in phenomena" [CP 1.418]. For these purposes, we can establish that a quality has a character—considered for its own sake—as a single but partial determination of what is capable of being felt. As such, it is an in eliminable part of our experience of any sort of phenomena.

On Peirce's account, the second category of the elements of the phenomena we observe comprises the actual facts about objects that react with other things or with us in some brute way [CP 1.418]. Qualities are vague and have the character of what is potential. They embody what could be or what is possible. Actual occurrences however are perfectly individual. Each occurrence happens at some particular time and place. With respect to our experience of something taking place before us at the present time, a brute fact happens here and now. Proverbially speaking, individual facts are called brutal because they resist our will. Those facts that persist in a more permanent sort of way are less purely individual, but the permanence of the fact persisting over time is a mark of its generality. We abstract in a precise way from such general features when attend to what is second in the brute actuality of the fact considered as an individual occurrence.

Peirce points out that "qualities are concerned in facts but they do not make up facts" [CP 1.419]. That is, every fact involves the attribution of some general qualities to some particular object. As such, facts have a dual nature in virtue of which what is potential as quality is brought into a relation to what is individual as actual object, but the hallmark of what is brute as fact is the object resisting any attempt on our part to bring about a change in the object.

On Peirce's account, the third category of elements of phenomena "consists of what we call laws when we contemplate them from the outside only, but which when we see both sides of the shield we call thoughts" [CP 1.420].  On this account, thoughts are not qualities because they can be produced and grow. What is more, thoughts must have some reasons, good or bad. In general, our experience of thought shows that they are the kinds of things that have the character they do in relation to other thoughts that serve as interpretations of the earlier thought--and these relations are some respects rational and are not merely brute relations. Or, to put the matter in different terms, thoughts are not brute facts because they have the character of what is general. If you are having a thought, it is the kind of thing that can be communicated or shared with another person. As general, thoughts do not refer solely to a limited range of actual things that happen to exist here and now. Rather, thoughts also refer to a wider range of possible things that, taken together, are neither discrete nor finite in this respect.

If we grant Peirce the claim that thoughts and laws are both examples of this third material category of experience, then a number of the same points that hold for thoughts hold also for laws. Laws are more than just a finite collection of facts in that they determine how a range of possible facts not yet realized may or must be in the future. The key point is that thoughts and laws concern both the realm of potential qualities and the realm of actual facts. As general, thoughts and laws govern the regularities of the qualities that can be found in facts when they are realized. As an element that has the character of what governs such relations between qualities and facts in a general way, this third element in experience is different in kind from the qualities and actual objects that it, in some sense, governs.

Other philosophers, ranging from Aristotle to Kant and Hegel, have offered competing accounts of the material categories in experience. Peirce argues that none of these philosophers have provided an analysis of what is really fundamental as a set of elements in experience. As such, Peirce pushes this line of inquiry further by carefully examining the elemental *formal* relations that are necessarily part of any experience—actual or possible.

Peirce says that our method must be to "observe how logic requires us to think and especially to reason," and to attribute to the conceptions of the elemental relations those characters which they "must have in order to answer the requirements of logic" [CP 1.444]. This, I believe, is key to understanding how the phenomenological method can be used to analyze the formal elements in experience. When Peirce makes the transition from the analysis of the material categories to the formal elements in experience, he needs a more refined method. In order to determine what kinds of formal elements—if any—are necessarily part of any possible experience or thought, we need to consider those aspects of experience that are essential for the possibility of using signs in a meaningful way and for making inferences from such signs.

The analysis of each of the formal elements that are necessarily part of any experience requires a clear separation between each of the elements. As such, he starts by telling us what each of the elements is not. The general idea he exploits is that each of the formal categories is *not* composed of the other elements, or the category in question would not be elemental. Having indicated what each is not, he then characterizes those aspects that are positively indicative of what each category is. As such, the analysis proceeds by focusing on the negative first and then turning to the positive characterization.

Peirce points out that, for those aspects of experience that are truly elemental, it is probably better to think of them as tones of thought than as conceptions [CP 1.355]. After all, the conceptions we employ are complex things that appear to be made up of many parts. So, what are these elemental tones of thought that are formal in character? Peirce explains that the simplest of these elements has the formal character of a monad. Upon analysis, we find that the simplest qualities in experience have the predominant character of what is monadic as a relation. As such, the experience of the color of red, the taste of what is sour, and the sound of a trumpet are, in their most elemental aspects, monadic as qualities. That is, if we separate characteristics such as the hue of the color, or the intensity of the taste, we find that one of these characteristics—such as the hue—considered in complete separation from all of the other characteristics, has its own monadic character.

The monadic character of such tones of experience is reflected in the signs we use to talk about such characteristics. As such, if we assert that "the fire engine is red," the predicate "is red" is a one-place relation considered from a logical point of view. When we move from a logical analysis of conceptions to a phenomenological analysis of ordinary experience, we see that the monadic element has the character of what is merely potential [CP 1.424]. When we engage in acts of precisive abstraction and attend to one element in experience and ignore the others, we must not even attend to the fact that we are aware of any determinate absence of other things. That, too, is something we must ignore when we attend to what is most characteristic of this elemental tone of thought. Peirce tells us that, in doing so, "we are to consider the total as a unit" [CP 1.424].

This, I believe, is a central point. Monadic elements in experience function as a unit insofar as we consider a given quality as a unitary whole. If we consider a given quality in experience as itself having further qualities—which many might very well have—then we are no longer attending to our experience of this tone of thought as a unitary whole. Setting to the side the question of whether every quality in experience is or is not composed of further characteristics—as a experience of color has its hue, chroma and intensity—we should note that he is asserting that each of these qualities of experience are when, considered as a whole, unitary in character. As such, Peirce's claim is that "quality is what presents itself in the **monadic** aspect" [CP 1.424].

This provides us with the first part of an answer to the question of how we should understand the basic feature of that number which functions as a unit. In our common use of the counting numbers the unit is the number one, and the same is true for the system of the natural numbers. Considered in its monadic aspect, the formal element of the monad in the conception of the unit should be considered a whole that takes what is *potential* in its character— and not what is actual--as what is most elemental. This feature of what Peirce is saying about the monadic character of the number one will be important later when we try to understand how the pairing of the monadic character with the number one fits with the account of the formal elements that are most prominent in the numbers two and three.

The second formal element of experience is the dyad. Unlike the monad, which is a formal element having the character of a unity that is considered without any regard to a structure or parts, the dyad is understood to have a structure that presents a variety of features. The characteristic material aspect of experience that has the formal character of the dyad is brute reaction. As we have seen, in the experience of brute reaction there is an experience of one actual thing that is resisting another actual thing. We experience such reaction when, for instance, we push to open a door that is ajar. Expecting the door to open, we exert a force and are surprised when that is met with an opposite force. In this balking of our will, we experience resistance to our efforts, and we understand such resistance in terms of a reaction between our efforts, as agent, and the door that is ajar, as the patient.

The formal element of the dyad is the underlying structural relation between distinct individual things, where one is resisting the other. In this way, the formal character of a dyadic relation is fundamentally different from the formal character of the monad. Whereas the monad has no structural relations or parts, the genuine dyad consists of a brute relation between two distinct objects. As such, the claim is that brute reaction is what presents itself in the **dyadic** aspect [CP 1.428-9].

Given the fact that the dyad can involve a number of different sorts of relations as part of its structure, we will need to take a closer look at the way Peirce classifies varieties of dyadic relations. In fact, the bulk of the discussion of the formal elements of experience in "The Logic of Mathematics" is devoted to a rich exploration and classification of the different sorts of dyadic relations might hold between two things. Peirce claims that the number two is intimately connected with this formal element of experience. As such, we will need to take a closer look at manner in which he classifies the different kinds of dyadic relations that might obtain between distinct subjects in our experience.

The third formal element of experience is the triad. We should remember, once again, that logic is our guide in this inquiry. Peirce points out that thinking of the triad as involving three subjects, just as the dyad involves two subjects, is to take an incomplete and somewhat misleading view of the matter. Logical analysis shows that the monad is embodied in the verb or the predicate. As such, we say of something that it "is red" or "has the monadic quality of redness." The dyad manifests itself in the subject of the proposition, such that we say "this thing" (e.g., a stop sign) has the property of being red. Just as the monad is correlative to the term signifying the predicate, the dyad is correlative to the sentence signifying a proposition that asserts something is actually the case. The triad, on the other hand, expresses the three-place relation that is found in the following expression: "the police officer gave the driver a ticket for running the stop sign." The relation of 'A *gives* B to C' expresses the logical character of the triad [CP 1.345, 8.331].

Let us ask: what is the formal element found in experience that is essential for the validity of reasoning about such matters? This is where the phenomenological part of the inquiry tries to answer a need of the philosopher who is developing a logical theory of the validity of reasoning. We can ask: was the ticket well deserved? That is, did the driver of the car run the stop sign, or did the officer fail to see that the car in fact came to a complete stop before proceeding? After all, it is possible that the officer's attention was diverted when the car came to a stop. The question of whether or not the act of giving the ticket is warranted or not is determined, in part, by the law that governs the case. In this way, we arrive at the idea that the triadic element is embodied in our experience that there is a *law* that governs the case, and in our experience that there is a *thought* involved in applying the law to see if the giving of the ticket was warranted by the actions of the driver.

It is incomplete to say that the formal element of the triad involves three subjects because the law and the thought are different in character from the manner in which we experience something as a simple subject. The law connects and binds the police officer, as one subject, and the driver of the car, as a second subject, to the ticket, as a third subject, in a manner that makes the action of giving a ticket more or less reasonable in a particular case. On the basis of a phenomenological analysis of this element of experience, we come to the conclusion that thought (or law when viewed from the outside only) is what presents itself in the **triadic** aspect [CP 1.471-2].

It is clear that the material category of *thought*, as it embodies the formal element of the triad, has considerable structure connecting the parts of the relation. As with the dyadic relation, what we need are some tools for more thoroughly analyzing the different kinds of structural and functional relations that hold between the parts of the triad. In the next section, we will begin the development of the needed tools. First, we will examine Peirce's logical definitions of the conceptions of what it is for something to be a *relative*, standing in a *relation* to other things, such that there is a *relationship* that is capable of bringing the parts together. Peirce is drawing our attention to the different parts of larger structures involving the elemental relations of the monad, dyad and triad. Having examined the logical definitions of these conceptions, we will then develop a simple graphical system that is meant to show how the figures Peirce uses to illustrate the character of the monad, dyad and triad can be combined to form a wide range of larger molecular structures.

**3. Relatives, relations and relationships: a simple diagrammatic system**

In the "The Logic of Relatives" (1897), Peirce sets the goal of clarifying the meanings of the words "relative", "relation" and "relationship" [CP 3.288-345]. Peirce applies the pragmatic method to these conceptions for the purpose of working his way up from a first grade of clarity to a second and then a third grade in the clearness with which we are able to grasp the meaning of these conceptions. Here are the three grades of clearness that Peirce characterizes in our apprehensions of the meanings of words.

1. The first consists in the connection of the word with familiar experience.
2. The second grade consists in the abstract definition, depending upon

an analysis of just what it is that makes the word applicable.
 3. The third grade of clearness consists in such a representation of the idea that fruitful reasoning can be made to turn upon it, and that it can be applied to the resolution of difficult practical problems. [CP 3.457]

In this essay, Peirce starts with the set of conceptions as they are drawn from familiar experience and our usage of the terms, and he provides nominal definitions that accord with familiar usage. Next, he provides the following abstract definitions of the key terms.

   a. A *relative,* then, may be defined as the equivalent of a word or phrase, which, either as it is (when I term it a *complete* relative), or else when the verb "is" is attached to it (and if it wants such attachment, I term it a *nominal* relative), becomes a sentence with some number of proper names left blank.
   b. A *relationship, or fundamentum relationis,* is a fact relative to a number of objects, considered apart from those objects, as if, after the statement of the fact, the designations of those objects had been erased.
   c. A *relation* is a relationship considered as something that may be said to be true of one of the objects, the others being separated from the relationship yet kept in view. Thus, for each relationship there are as many relations as there are blanks [CP 3.466].

With these logical definitions in hand, let us apply them in the development of a simple graphical system that is based—at least in spirit—on Peirce's account of the elemental relations of the monad, dyad and triad and his explanations of how those elementary relations can be formed into more complex molecules. The ground rules for understanding how the elements can be combined are based, in part, on Peirce's criticisms of Alfred Bray Kempe's dyadic model for analyzing mathematical forms of expression [Kempe 1886, 1887, 1897] [MS 708-9]. We won't review those criticisms here because it would take us too far afield from our main interests, which is to show how this particular diagrammatic model works. The main idea that animates this model is to preserve the distinctive features of the elemental monadic, dyadic and triadic relations as they are characterized in the phenomenological theory.

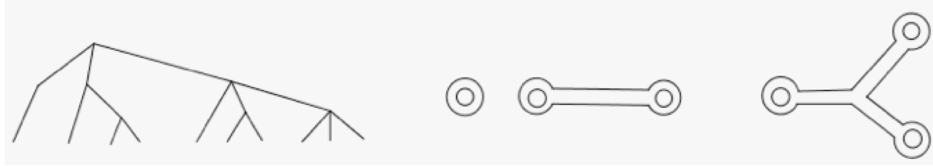

> **Figure 1.** Peirce uses two different kinds of diagrams to illustrate the three formal elements of the monad, triad and dyad: branching tree diagrams and roadways with turnabouts.

In his later works on the phenomenological categories, Peirce introduces two kinds of diagrams that are illustrated in figure 1 [EP 162, CP 1.363]. In the diagram on the left side of the figure, there is a branching tree structure that shows a monadic relation branching from a triad and various degenerate and genuine dyadic and triadic relations branching further. On the right side, there are a set of roadway diagrams with a turnabout: first, a simple turnabout; second, a road with two turnabouts; third, a branching road with three turnabouts. Each of the roadway diagrams represents a fully saturated medad having no free ends. The roundabout itself is a simple medad with no parts. The roadway with two turnabouts is composed of a dyad and two monads, and the roadway with three turnabouts is composed of a triad and three monads.

Our goal is to develop a simple system that is modeled on a graph theoretical approach in mathematics [Biggs, 1976] [Beineke 1997] [Grattan-Guinness 2002]. Peirce thinks of graph theory as being a part of the larger field of topology, and he suggests that these areas of mathematics are particularly important for the study of true continuity [Havenel 2008]. Simple tree branching diagrams were used by Peirce to explore questions in mathematics and formal logic, so let us describe a system of rules for creating and manipulating diagrams that is motivated, at least in spirit, by Peirce's own inquiries using similar diagrammatic structures [Anellis 1990] [Anellis 2016].

Let us start with three basic elements: a monad, a dyad and a triad.

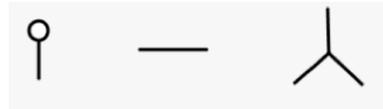

**Figure 2**. From left to right, there is a diagrammatic representation of a monad with one free end, a dyad with two free ends, and a triad with three.

These are illustrated in figure 2. The monad is represented as a circle with a single line. The circle represents a bounded area that cannot be joined to other things, and the tail end of the line segment represents one open end. The dyad is a line segment with two tail ends that represent two open ends. The triad is a branching figure that represents a connection between three line segments (i.e., a bifurcation)—each of which has an open end. Connecting these three elements at their open ends will compose all other figures.

The basic rules for forming and breaking connections between these elements are as follows:

(1) An element may be connected to any other element at their respective open ends. When a connection is formed, let us designate the bond with a diamond. Any time two or more elemental relations are joined together, the resulting object is called a molecule

(2) Each bonding site is fully saturated when one open end has been connected to another open end. It is not permissible to connect more than two elements at their open ends.

(3) A connection between two elemental relations may be removed by erasing the bond at a saturated site so that the parts are disconnected.

Let us start with some examples of the simplest kinds of connections that can be formed between elements. In figure 3, we have three different kinds of molecules.

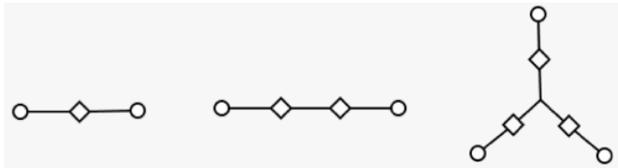

**Figure 3**. This diagram shows three permissible types of molecules that are composed of formal elements. On the left, two monads have been joined. In the middle, there is a dyad with two monads. On the right, there is a triad with three monads. The sites where the open ends of the elements are bonded are represented with a diamond. Each molecule is fully saturated and functions as a medad.

The diagram on the left of figure 3 is a relationship between two monads. The figure in the middle is a relationship between two monads and a dyad, and the figure on the right is a relationship between three monads and a triad.

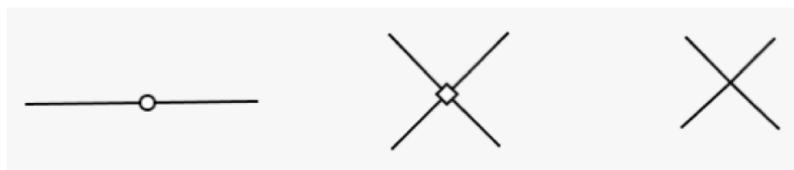

**Figure 4**.Three examples of impermissible forms of combination of formal elements.

In figure 4, we have three examples of connections that are not permissible in this graphical system. The first is a monad that has more than one line segment radiating from the circular portion. The second is a particular bonding site that connects more than two open ends. The third is a branching figure with more than three line segments radiating from the furcation.

    Let us introduce a bookkeeping idea that will help to keep things clearer as we put the graphical system to work in studying the different classes of structures that can be formed by combining elements and molecules. We can mark a division between relations that are internal and those that are external using a dashed circle.

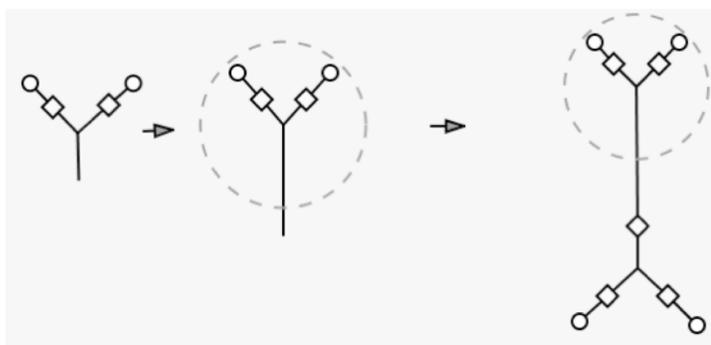

**Figure 5**. Dashed circles are used to mark a combination of elements. The inner part of the circle shows the relations that are considered to be internal. In the middle, a triad is marked with a dashed circle to show one free end, which is then bonded to another triad of the same type.

In figure 5, there is a molecule on the left side of the diagram that consists of a triad connected to two monads with one open bonding site. Before any more connections are made between this molecule and other elements, let us draw a dashed circle around the element so that only the unsaturated bond extends beyond the circle. The purpose of designating the difference between internal and external relations with the dashed circle is simply to keep track of the parts and connections at each step in the process of adding and removing pieces from the structure. It is clear that, for some purposes, the area enclosed by the dashed circle will function as a monad. The portion of the molecule that is enclosed by the dashed circle does have internal complexity but, for the purposes of adding additional elements, there is only one open end.

It is worth noting that the triad serves the unique function of enabling us to connect saturated molecules to other molecules by a process of insertion of the triad into the structure to create new open ends. In effect, it enables us to form connections between more complex relations so as to create new unsaturated bonding sites. Neither of the other two elements are able to serve this critical function.

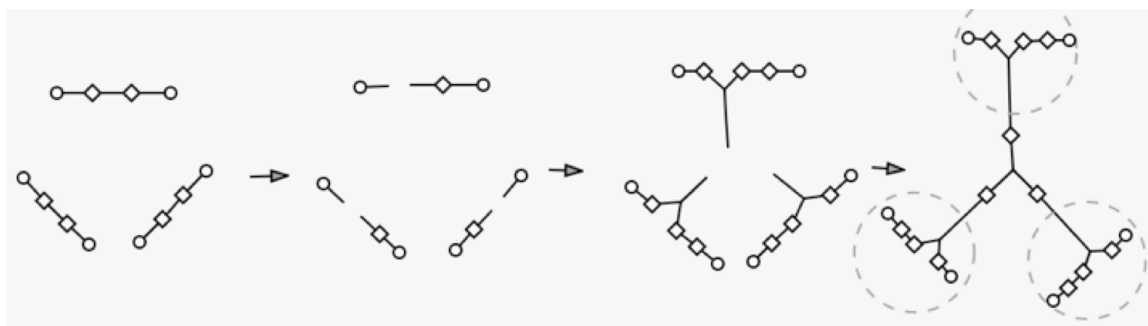

**Figure 6**. On the left of the diagram there are three saturated dyads. The process of combination involves the breaking of the bonds within each of the dyads and the insertion of a triadic element. As a result each of the three new structures has one free end, which are then joined to a triadic relationship to form a new molecule that has the structure of a triad of dyads.

In figure 6, a process of combination starts with three saturated dyadic molecules. In order to form connections between these molecules, it is necessary to erase one of the saturated bonding sites in each of the relations. Now that there are open bonding sites, a triad can be inserted between the disconnected pieces of each of the dyads. In this manner, each of the dyads can be connected to the others. The result is a triad of dyads.

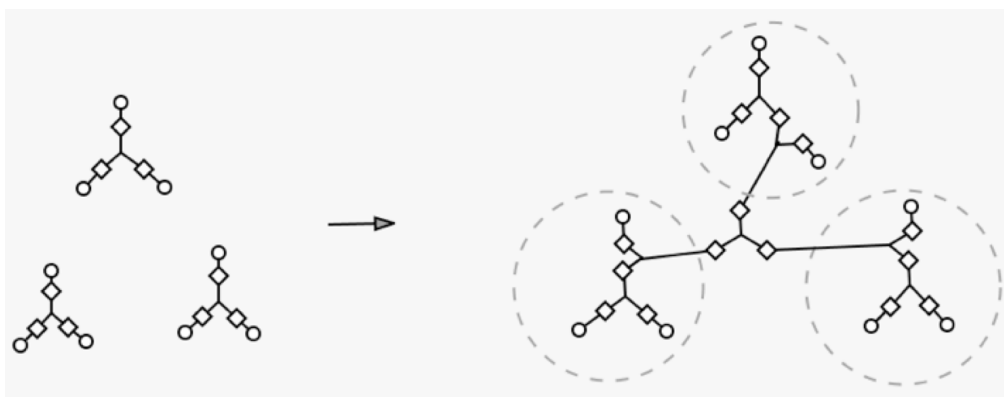

**Figure 7**. This diagram illustrates the process by which three saturated triadic molecules are joined by a triadic relationship to form a molecule having the structure of a triad of triads.

The same point holds for other kinds of relations, including three saturated triads. By the same process, we can form a triadic relationship between three triads (see figure 7).

Let us note that, after a molecule has been formed, certain parts of the diagram can be erased with only a small loss of information. In particular, the diamonds that represent the saturated bonding sites can be erased along any continuous line leaving the connection intact. For any continuous line in a molecule, there is at least one bonding site between monadic terminal ends or bifurcations. The only information that is lost is the number of dyadic relations that were connected in a linear row. Given the fact that these connections between dyads can be inserted or erased without altering the basic structure of the diagram, this often will turn out to be an immaterial loss of information. As such, I will often drop the use of the diamond in those diagrams where it is apparent how the monads, dyads and triads have been connected in previous stages of the composition of the relations for the sake of highlighting the new bonds that have been formed.[3]

Before going further, I should note that the diagrams of *triads of dyads* and *triads of triads* illustrated in figures 6 and 7 will be essential for interpreting Peirce's explanations of how monadic, dyadic and triadic relations are involved in the numbers one, two and three. In particular, two different kinds of triadic structures will be central in our interpretation of what Peirce is doing when he provides an analysis of the formal categories involved in the number three.

## 4. Applying the graphical system to the classification of relations

---

[3] One reason to retain all of the parts of the figures—including the diamonds that represent the saturated bonds—is that it makes it possible to study more closely the combinatorial possibilities of the system.

The graphical system developed above can be used to analyze the various ways that monadic, dyadic and triadic elements can be combined to form larger molecular structures. Let us use this system to provide an explication of what Peirce says about the classification of dyadic relations. Having done that, we will then turn to the division between the main classes of triadic relations.

The dashed line in figure 8 (below) signifies a division between those dyadic relations, found in the top half of the figure, whose classification is based on the kinds of *subjects* that stand in the relations [CP 1.454]. The classes below the dashed line are classified not on the basis of the kinds of subject*s* that stand in the dyadic relations but on the kinds of dyadic *relations* that hold between the subjects [CP 1.456].

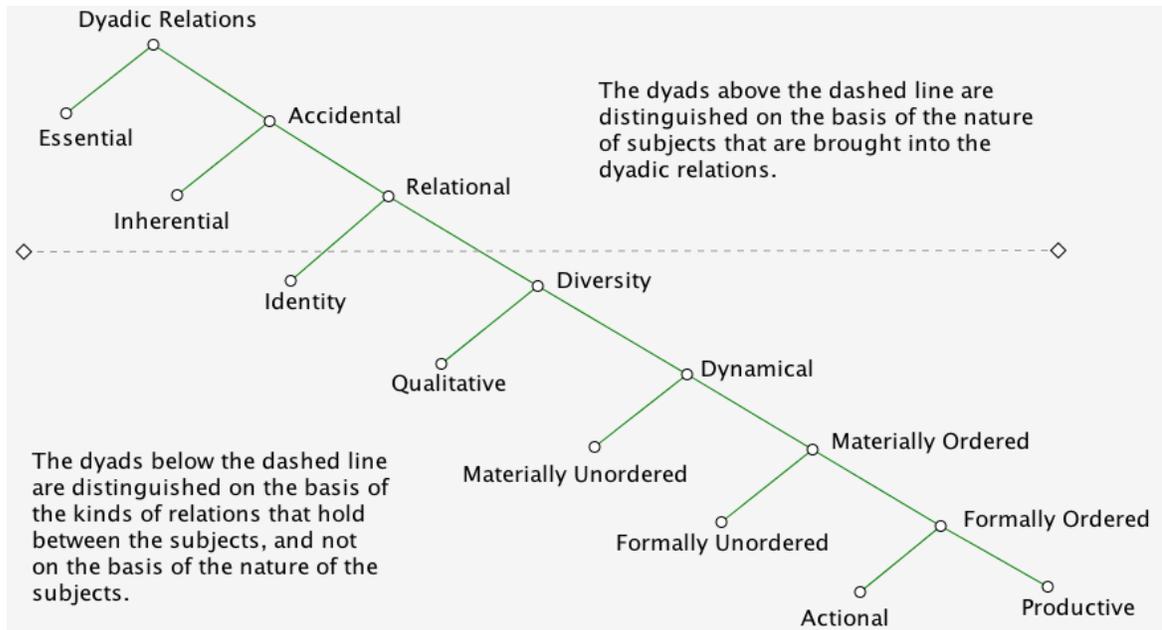

**Figure 8**. This diagram shows the tree branching structure of the different classes of dyadic relations. The character of the subjects in the dyadic relations defines the classes above dashed line, and the character of the dyadic relationship that connects the subjects defines those below the dashed line. At each branching point, a more degenerate class is on the left branch and the more genuine dyadic relation is on the right branch. Only the more genuine class is then divided further into additional classes.

Using the terminology Peirce employed in the discussion of the formal elements of the dyad, different dyadic relations can be classified based on the structural relations that hold between the parts that make up the subjects of the dyads, and they also can be classified based on the kinds of structural relations that hold between the subjects. As such, the divisions between the classes above the dashed line are based on the *internal* structure of the subjects that are brought together by the dyadic relationship, and the divisions between classes below the line are based on the character of the *external* relational structure that connects the subjects.

As we have seen, a monadic relation is embodied in a simple quality, such as the character of "redness", considered as an abstract potentiality. In a certain respect, the character of redness is considered as an abstraction because that is the manner in which we have arrived at this tone of thought as an element of experience—i.e., we have arrived at it by a process of precise abstraction. Considered as a simple element, we consider it as a unit without any parts or structure.

The classification of dyadic relations can be arranged on the basis of a pattern of branching relations. The branches to the left at each level in figure 8 represent relatively degenerate forms of dyadic relations. The branches that continue down to the right, relatively speaking are, more genuine forms of dyadic relations. The more degenerate relations represent dyads that are, in some respects, more monadic in character. The more genuine relations represent dyads that are, in respects more dyadic in character. Let us examine the diagrams in order to see how the classificatory system based on degrees of degeneracy and genuineness works.

The most degenerate kind of dyad is one for which the loose ends of the dyadic relationship have been saturated with two subjects that are simple monadic qualities. As such, the dyadic relationship that holds between of the quality of scarlet and the quality of red has the character of the one being contained within the other as an essential dyad. This type of degenerate dyadic relation is classified as an essential dyad because the relation of being "contained within" is a relationship that holds between the two qualities that are differentiated solely based on the essential character of each considered as simple monad [CP 1.455]. As simple monads, the qualities are being considered as subjects that have the character of abstract potentialities.

Accidental dyads have at least one subject that is something more than a mere monad. Inherential dyads, which are the more degenerate class, have one subject that is a mere quality, while the other subject is a dyadic relation holding between two subjects [CP 1.456]. Relational dyads, which are the more genuine class, have two subjects each of which consists of one or another type of dyadic relationship (i.e., essential, inherential or relational), each of which is illustrated in the diagram.

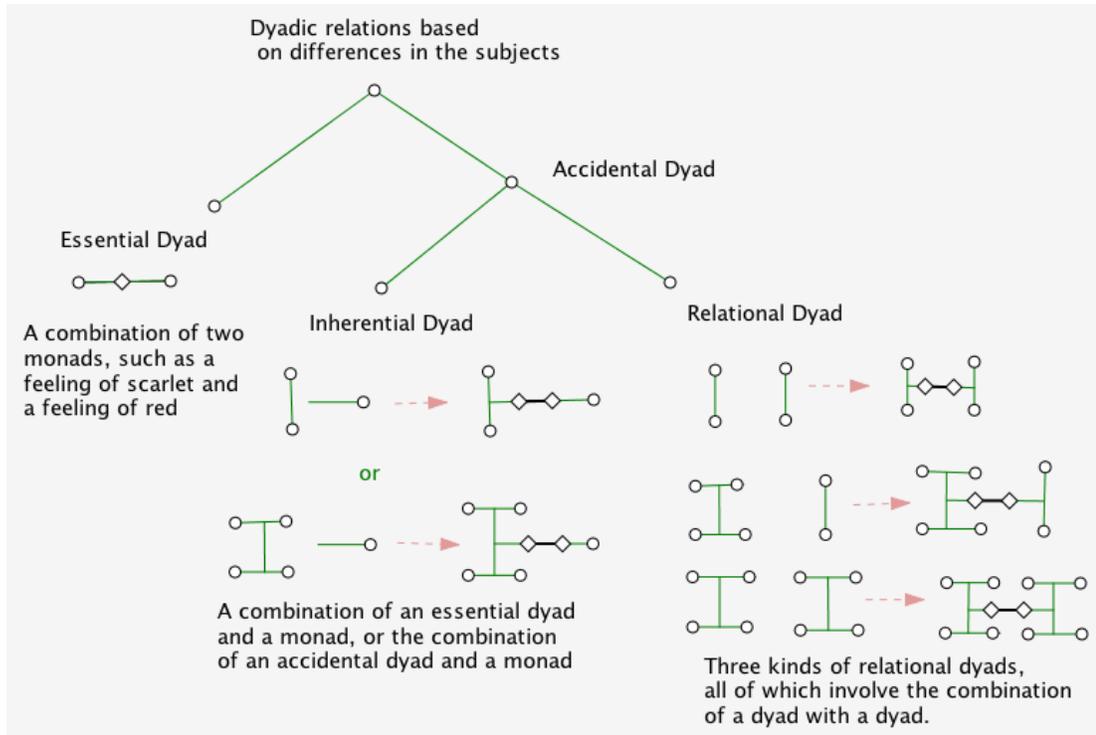

**Figure 9**. This diagram illustrates the process of creating the different kinds of dyads that are characterized by the nature of the subjects in each of the three classes of dyadic relations.

Having described the classes that are distinguished based on the nature of the subjects that are connected by the dyadic relationship; let us now consider the classes below the dashed line in figure 9. All of these classes are distinguished based on the nature of the relationship that connects the subjects, and all of the subjects are relational dyads of one kind or another. The most degenerate class of relational dyad is one of identity between the subjects, such that the two subjects are really one and the same individual thing. All other dyads are classes of dyads of diversity such that the two subjects are not the same individuals. The dyads of diversity are divided into qualitative dyads, which are the more degenerate class, and dynamical dyads, which are the more genuine class. Qualitative dyads involve two inherential dyads that are combined such that one subject possesses a quality that is different in character from the quality of the other subject [CP 1.464]. Dynamical dyads of diversity, on the other hand, are more genuine because they involve an object having a quality that is brought into a relationship with another object having a quality such that one of the objects is dynamically giving rise to a change in one or more of the characteristics of the other object [CP 1.465].

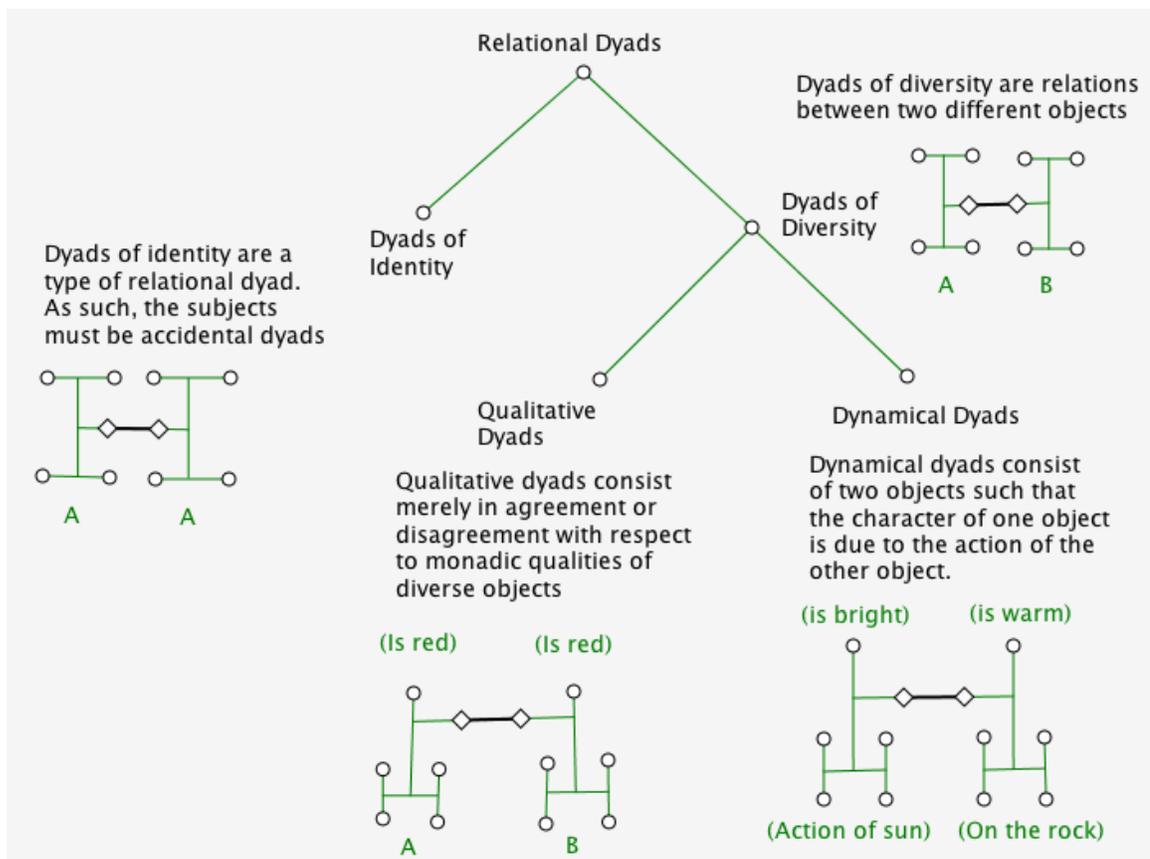

**Figure 10**. This diagram illustrates the structure of the different classes of relational dyads.

Dynamical dyads are further classified according to the kinds of material or formal order that holds between the subjects, and finally according to the kind of ordered action that one subject exerts on the other subject [CP 1.467]. Looking at the classes of dynamical dyads, it would appear on its surface that the order found in the dyadic relations that are typical of materially and formally ordered dyads are, in some way, important to gaining a better understanding of how the dyad is connected with the number two. What is more, the order would seem to be important to better understanding the role of the conception of two in the larger set of hypotheses that lie at the bases of discrete systems of mathematics.

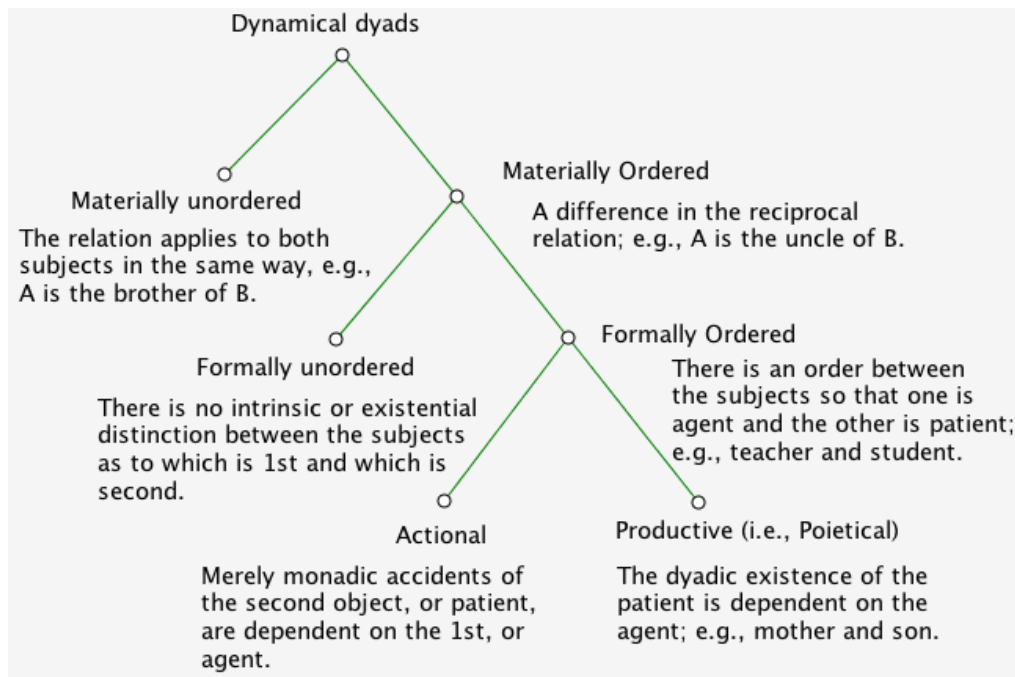

**Figure 11.** The branching structure in this diagram highlights the more complex structures of the different classes of dynamical dyads. The different kinds of ordered relations may have special significance for explicating the character of sequences of collections of objects.

Unfortunately, Peirce does not explain in any detail what significance such ordered relations have for our understanding of the conception of two or the hypotheses that we are trying to explicate in this essay.

What does seem to matter for the sake of understanding the connection between the number two and the dyadic relation is that *correspondence* is a relation that holds between actually existing individuals. Only the experience of such actually existing individuals has the character of what is discrete, so this type of experience is the source from which we draw our conception of the number two. With that much said, let us proceed to the analysis of triadic relations and then return to the question of the role different kinds of ordered relations might have in supplying answers to the questions that were posed at the start of the essay.

## 6. The Classification of relations

On Peirce's account, there are three major divisions in the basic classes of triadic relations [CP 1.473]. The most degenerate class consists of those triads in which each of the subjects is a simple monadic quality. As such, the most degenerate class of triadic relations has the same kind of subjects as the simplest class of dyadic relations. That is, the triadic relations are essential in character. Peirce offers the example of scarlet being red which in turn is a color. The relations between scarlet and red is that of being contained within, just as the relation of red and color is that of being contained within.

The next class of triadic relations consists of those in which each of the subjects is a dyadic relation. Given the classification of the different kinds of dyadic relations that are part of experience, it is clear that there are a whole host of different kinds of combinations in which one subject might be a qualitative dyad of diversity, and another might be a dyad of identity, and yet another might be a formally ordered dynamical dyad.

**Figure 12.** This diagram illustrates the three main classes of triadic relations. The two degenerate

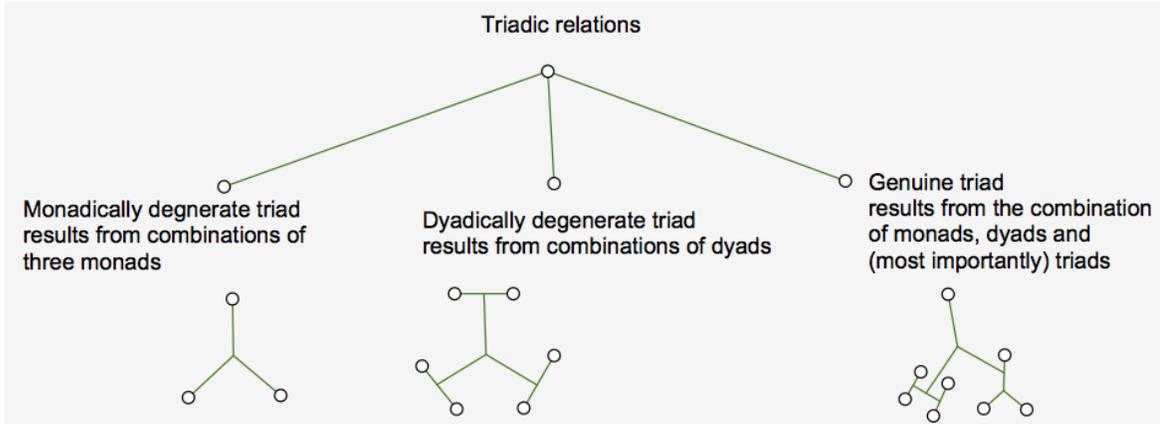

classes seen on the left and middle include those that are composed of three monads joined by a triadic relationship, those composed of three dyads joined by a triadic relationship. The classes of genuine triadic relations on the right have at least one triadic relation as the third correlate in the relation.

It would appear that different kinds possible combination of dyads that can be the first, second or third correlate would depend upon the character of the triadic relation that brings them together. Having said that, we can probably move forward without digging into the details about how such dyads might be combined in this degenerate class because it is the genuine triadic relations that matter most for what he says about number.

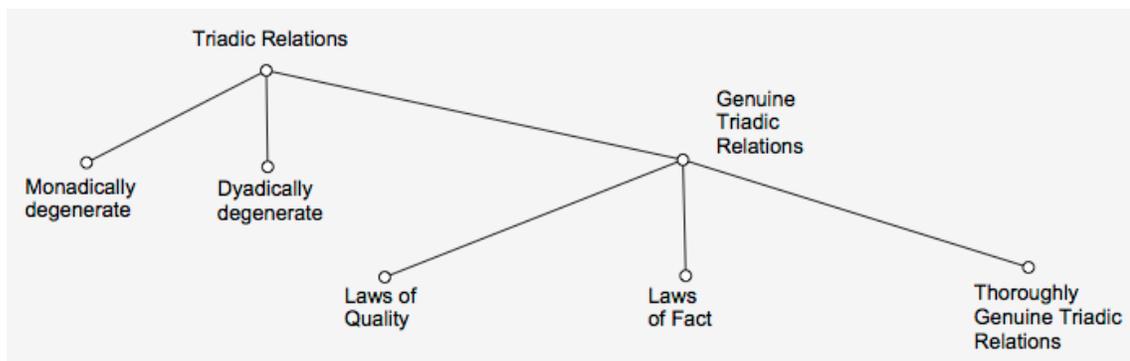

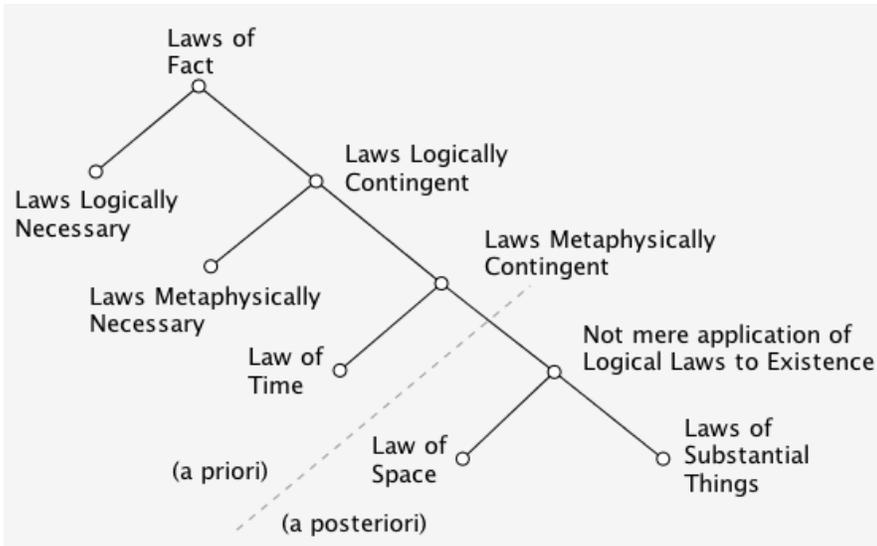

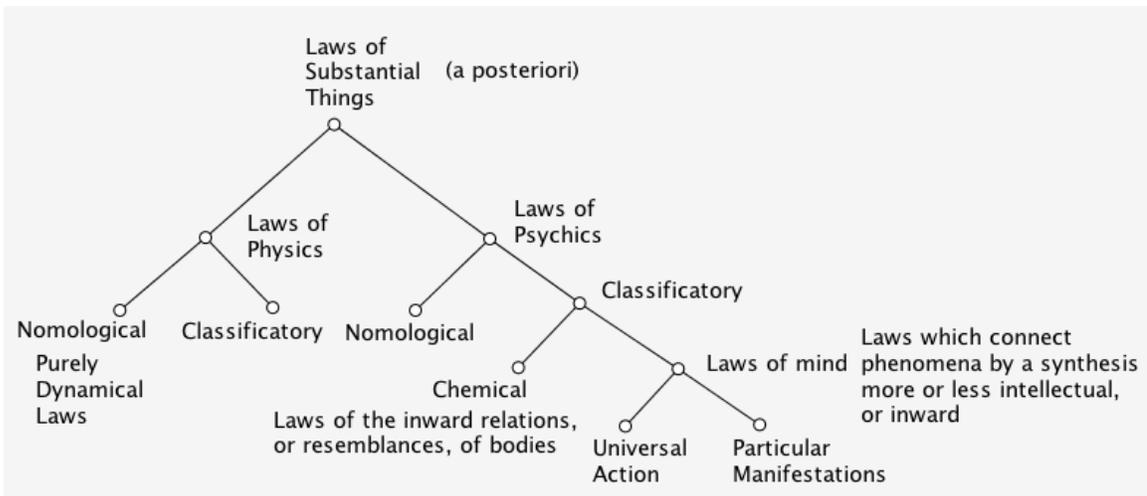

It would be wrong, Peirce says, to "define two as the sum of one and one" [CP 1.476]. The reason it would be wrong is that the main idea involved in the number two is that of other. This idea of *other* is prominent, for example in a dynamical relation of diversity between two existing objects. The ideas that are prominent in the number three are the ideas of *third*, *medium* or *uniter*. While the genuine triad involves these three ideas, the defining characteristic of the genuine triad is that of a "third not resoluble into a formless aggregation" [CP 1.477] As such, the genuine triad involves the idea of something that is more than what can result from a dyadic addition of one and one. The idea of "something more" involves the notion of "every possible something" and hence the idea of generality. As such, the fundamental difference between the number two and the number three is that the latter involves this generality that applies to an unlimited range of possibilities.

**6. Collections, counting and the precepts for number systems**

At this point, let us take a step back from the phenomenological analysis and summarize what we learned thus far. If Peirce's arguments are on track, then we have good reason to believe that the formal categories of the monad, dyad and triad are necessary elements in all possible experience. In particular, these features are essential components in the observations that we might make of any sort of phenomenon that might call out for explanation. As such, these formal elements are entirely necessary for reasoning from observations to conclusions, including observations of the formal features of idealized diagrams. Consequently, they are the formal elements from which the basic conceptions of one, two and three are comprised.

Phenomenological analysis of the numbers one, two and three shows us the

following:

1) The number one is correlative to the monad as a formal element insofar as it is thought of as a *unit* that pertains to what is *potential*. That is, the number one is not being considered as correlative to an actual subject that is individual and discrete. Rather it is correlative as a unit to an unlimited range of what is possible.

2) The number two is correlative to the dyad as a formal element insofar as we are considering richer types of dyads, such as what are found in the experience of subjects that are diverse, relational and stand in dynamical relations. The reason is that the number two involves the notion of one discrete individual thing standing in a relation of correspondence to another discrete individual thing.

3) The number three is correlative to the triad as a formal element in experience insofar as genuine triads involve the notion of a general thought or law. The reason is that the number three involves the combination of different subjects under such a general law or thought. What is general covers an unlimited number of possibilities regardless of whether it is considered to have the characteristics of a thought or a law.

Let us apply these ideas to a specific set of claims Peirce makes in "The Logic of Mathematics" about the role of these formal elements in our understanding of what it is for something to be a member of a set. Peirce starts this line of inquiry by asking about the *function* of the units of a set in the *constitution* of that set. This is how he articulates what the units do. The first point he makes is that, in logic, "a set cannot generally be adequately represented by a diagram of a promiscuous collection of dots" [CP 1.446]. The reason is that a collection dots can be arranged in many ways. Some collections might form a disorganized cluster; others might be arranged in a linear row, while others might be arranged in the shape of a square. Peirce points out that the last kind of arrangement might be fitting for thinking about the character of a determinant.

In virtue of what property does a collection of a dots comprise a set? Peirce's initial point is that being a set is a not a property of the collection of the dots considered in themselves. For any given cluster of dots, a set might be formed from one of the members, or some larger portion of the dots in the cluster, or from all, or even from none. What makes a set the kind thing that we can reason about in mathematical inquiry is that the *form of connection* belongs to the set itself, and it does not belong to the units (e.g., the dots) that are taken to be members of the set.

Logical reasoning about the set requires that we consider this form of connection that makes the set the kind of thing that it is. As Peirce is keen to point out, all reasoning is formal in character. As such, any inference that is sound concerning a thing or character is sound in regard to any other thing or character, so long as the form of connection of the one inference is strictly analogous to that of the other inference. What follows from this general point about the character of valid reasoning is that, for the purposes of drawing logical inferences, all that has to be represented is the characters of the sets themselves. That is, "the units need exhibit nothing except what is requisite to the exhibition of the characters belonging to sets" [CP 1.466].

Peirce encourages us to ask: "What then, is the use of the units, at all?" Given the fact that we strip of all of material qualities of the object away by precise abstraction when we consider an object to be a member of a set, what does the conception of the unit contribute to the representation of the distinctive characteristics of a set? For the purposes of inquiry in metaphysics, we can focus our attention on the mode of connection of the parts that form a collection, and we can express that connection in abstract terms, without making any reference to the units that make up the collection. After all, in metaphysics we might ask what makes a collection the kind of thing it is, really? Is the collection really something that is more than just the sum of its parts, or is it not?

In mathematical inquiry, on the other hand, we have a different purpose in mind. Our aim in making a representation of a set is to connect it to the representation of another set. In order to make that connection, we must consider the units in virtue of which a thing can be a member of the different sets that are brought into this relation of attachment [Moore 2009, 2010]. As Peirce points out, "in order to show how the total set is composed of those two sets, it is necessary to take account of the identities of their common units." But, as we have seen in our examination of dyadic relations, *identity* is a special kind of relation. As a form of relation, identity cannot be implied by a general description of those things that are taken to be identical.

The descriptions of the sets considered in themselves apart from the individual things that are often taken to be members is a description of a thing that has the character of a general. This is clearly seen in the fact that we can form an idea of the null set: e.g., $\emptyset$, or a set with no members $\{\}$. What is more, we can form a set that has as its members the null set and the set of the null set: $\{\emptyset, \{\emptyset\}\}$.[4] The conclusion Peirce draws is that the only purpose in specifying the units in a given representation of a set is in order to establish that each unit of that set may signify its identity or lack of identity with an individual that belongs to another set. Or, to put the point more succinctly, the only function of the units of a given set is to establish *possible* identities with the units of other sets. In this way, we get a clearer understanding of why it is important to note that the number one is correlative to the unit considered as a potential sort of thing and not an actual thing that is considered to be individual and discrete.

---

[4] The operation of union in virtue of which such a set is formed can help us understand the relationship between the number two considered cardinally and the number two considered ordinarily.

Considering the two examples offered above the null set and the set that has as its members the null set and the set of the null set, we see that the latter set has one member that is identical with the sole member of the other set. In this case, the kind of potentiality in the unit that is important for being able to reason mathematically about sets is the possible identity of the member of one set with a unit of a different set. As we have seen, the identity of discrete individuals is a classified as a degenerate form of a relational dyad [CP 1.466]. As a class of relational dyads, it requires only two subjects. If we want to go further and determine whether a set that consists of three objects are identical, we only need to consider each of the three pairs of objects and ask if any of the pairs are identical. As such, the only kind of identity that matters for the sake of reasoning about the possible identity or diversity of the members of a set is degenerate relational dyadic identity.[5]

Having considered what Peirce says in "The Logic of Mathematics" about the different kinds of formal relations that are fundamental for reasoning about sets and their units, let us turn to his discussion of number in the *Elements of Mathematics* and a related manuscript.

---

[5] This account of unity and identity does not appear to be subject to any of the objections that Frege raises against alternate accounts drawn from the history of mathematics and philosophy. Furthermore, Peirce's explanation of the character of these conceptions seems considerably richer than the proposals made by Frege insofar as Peirce offers better explanations of how these conceptions are drawn from the formal elements found in common experience [Frege 1980, 24-95].

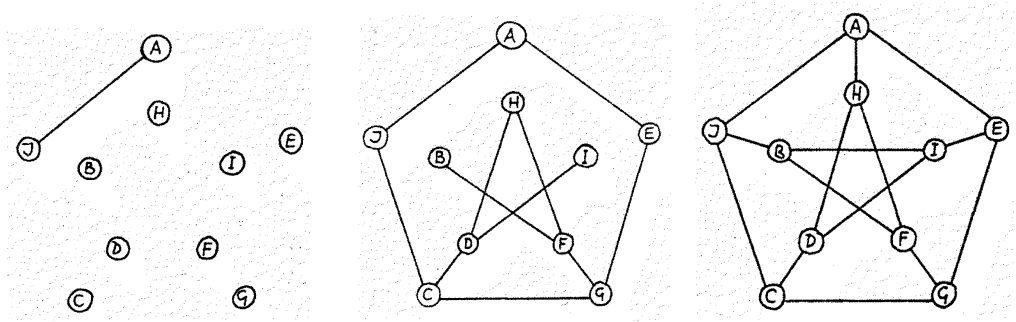

Figure 13

The example he offers is drawn from Kempe's work on the nature of mathematical form. Kempe's arguments about the nature of things such as sets and numbers draw heavily on a range of diagrams involves dots and edges, and the example is explored in a manuscript in which he responds directly to Kempe's analysis of these graphical relations. The diagram on the left in figure 13 shows a promiscuous collection of ten dots, with a line connecting dots A and J. Each dot, one might suppose, represents a monadic relation, and the edge connecting two spots would seem to represent a dyadic relation. Peirce argues that this is not a correct analysis of the predominant character of the parts of the diagram [CP 5.82-7, MS 708]. This is made obvious when we consider the way the diagram on the left is further developed as more edges are added to make the middle diagram. Eventually, a fuller pattern emerges, and each spot is connected to three edges. We can see that dots, as they are being used in this formal system, connect several things and do not have the character of a monad. Even in the first diagram on the left of the figure, the dots have the potential to connect several edges together—even when there are no edges radiating from a given dots. As such, the dots do not have the character of the unit in this system. Rather, they have the character of the triad. Each edge connects two dots and not more. As such, the edges are dyadic in character. Peirce argues that one thing in the system that does have the character of the monad is the blank sheet on which the spots are being scribed.

Let us clarify the claims Peirce makes in this argument by drawing a comparison between what has become the standard approach to representing these kinds of relations in graph theory to the manner we have been representing the elemental formal relations of the monad, dyad and triad in the diagrammatic system developed above. In order to make this comparison, we will need to formalize the graphical system and make the first steps in treating it as a pure system of mathematical relations. Let us call this formalized graphical system "MDT" graph theory (i.e., our graph theory of the elemental "Monad, Dyad, Triad") and note that, in developing the system, we are guided by the *aim* of learning to experiment with the molecular structures that can be formed by processes of composition, iteration and erasure of bonds between monadic, dyadic and triadic relations.

Towards this end, let us illustrate how the steps of constructing diagrams in standard graph theory can be represented in basic set theory. Then, we will illustrate how the steps of constructing diagrams in MDT graph theory can be represented in basic set theory. So, for starters, consider the process of connecting simple collections of dots and edges in linear chains.

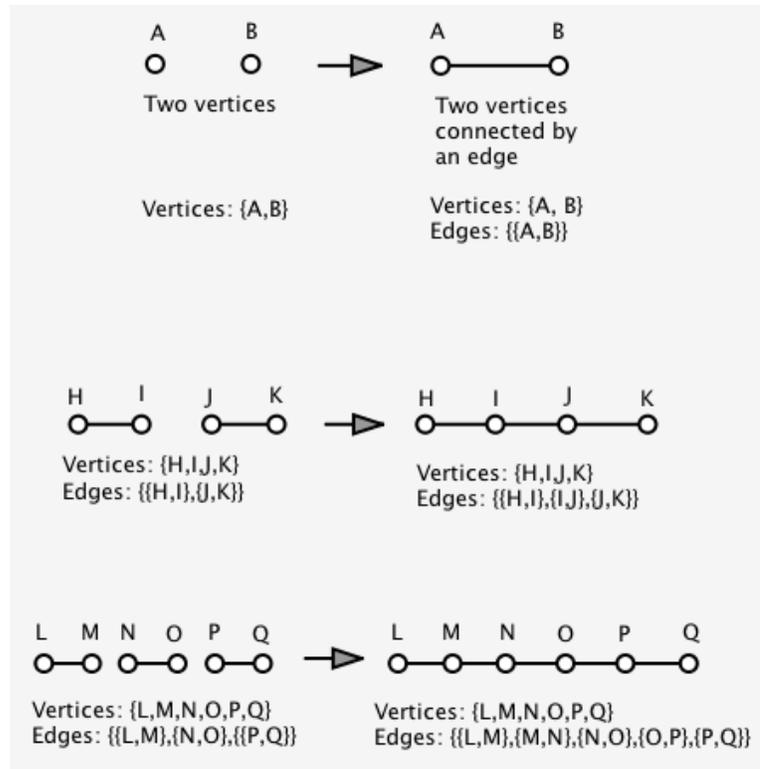

**Figure 15**. Connecting linear chains of edges and vertices in standard graph theory.

In figure 15, the diagram at the top illustrates how the process of connecting of dots with a single edge is represented in standard graph theory. The diagrams in the middle and bottom of the figure illustrate how the process of forming longer linear chain can be represented in basic set theory. In each case, only two types of sets are needed to represent the relations between the spots and edges: one type of set consists of vertices as elements, and the other type of set of consists of edges, where the edges are represented by a set of the sets of vertices at each end of a given edge.

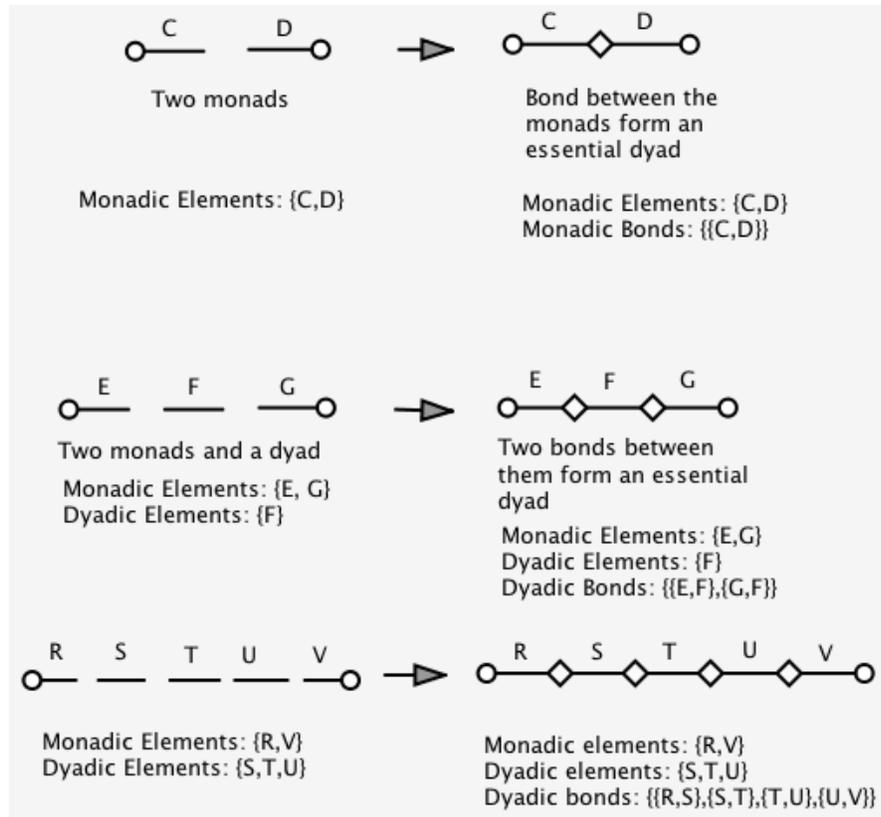

**Figure 16**. Connect linear chains of monads and dyads in MDT graph theory.

In figure 16, diagrams of linear chains consisting of bonds between monadic and dyadic relations are being formed in MDT graph theory. Note that the representation of the monadic relations and the bond that is formed between them is represented in terms of the correlative sets in a manner that is similar to the way the connection between one edge and two vertices in standard graph theory was represented in terms of the correlative sets (see figure 15 above). This, I suspect, is part of the reason many people, including Kempe, might be lead to assume there is no real difference between the two systems and, as such, there is no need to move beyond standard graph theory for the purposes of analyzing the relations between spots and lines in the diagram we have considered earlier (see figure 13 above.)

The same is not true when dyadic relations are represented in the diagrams in the middle and bottom of figure 16. In these two diagrams, it appears to be necessary to create a third and a fourth type of set in order to represent the dyadic relations and the bonds between them. While it might appear that two more types of sets are needed for dyadic relations and the bonds that hold between them in addition to those that are needed to represent the set of monadic elements and the bonds that hold between them, this may be misleading. The reason, as Peirce explains, is that a bond between two monads is, for all intents and purposes, an essential dyad.[6]

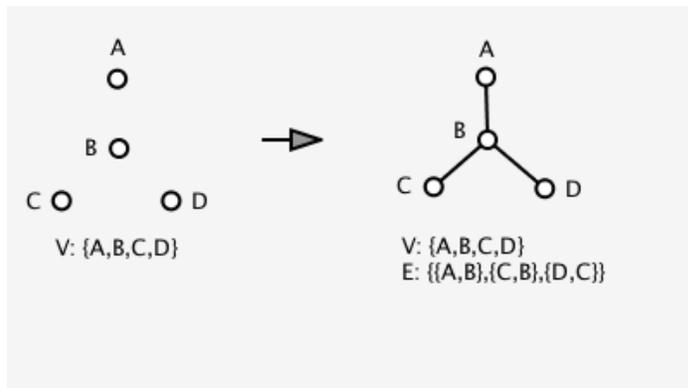

**Figure 17**. Forming triadic structures with vertices and edges in standard graph theory.

---

[6] For our purposes, a bond between two monads can be treated as involving a dyadic relation having the character "the hue of scarlet is contained within the hue of red." This kind of relation is made clearer by thinking about the relations between regions on a color chart where we can clearly see that any possible shade of scarlet should be found within the possible shades of red. Even if we can't draw determinate lines around the borders of the respective hues in terms of the relations of similarity, it is clear that the scarlet region on the chart is entirely contained within the red region on the chart. Having said that, it may be a confusion to think of hues of color—as abstract properties having the character of monadic relations—to be members of sets because the relations between such monadic qualities is not one of correspondence between discrete individuals—despite the insistence of Humean minded philosophers that all impressions of sense are actual, discrete and simple. This monadic formal feature in the experience of qualities of experience may be a source of confusion for those, like Nelson Goodman, who seek to build a phenomenological account of the basic structural relations that hold between the elements of experience using set theory as a formal model [Goodman 1966].

The diagrams in figure 17, illustrate the process of forming what look like branching triadic structures in standard graph theory. As with the case of forming linear chains, the processes of forming these structures can be adequately represented in simple set theory using only sets of edges and vertices. Consequently, there does not appear to be any significant difference between the way linear chains and the way more complex branching structures are represented in standard graph theory.

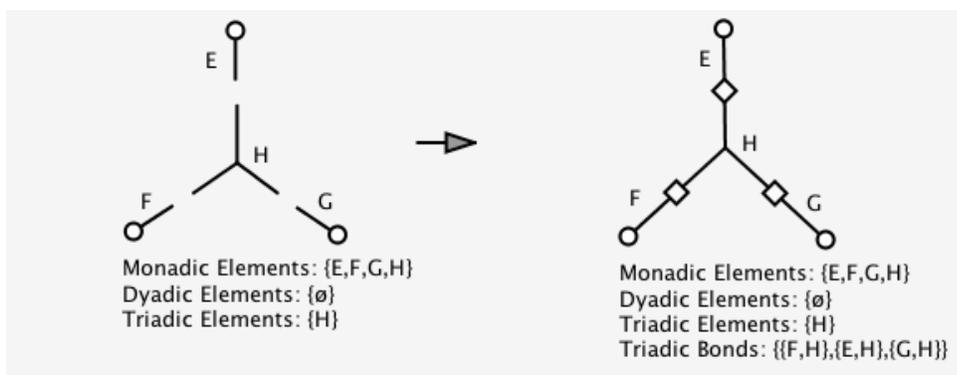

**Figure 18**. Forming bonds between monads and a triad in MDT graph theory.

The diagrams in figure 17 illustrate the process of forming triadic structures in the MDT graph theory. Note that four types of classes seem to be needed to represent the structure of the triadic molecule that has been formed as a result of the process. At this point, we begin to get a sense that there are significant differences between standard graph theory and MDT graph theory. The immediately apparent differences are the following. The diagrammatic system in standard graph theory seems to have two kinds of elements: vertices and edges. The diagrammatic system in MDT graph theory seems to have four kinds of elements: monads, dyads, triads and the bonds that join the free ends of these relations. These simple differences are not minor. The inclusion of a specific sign for the monad in MDT graph theory gives the system the capacity to represent a relation that is terminal in character. What is more, the inclusion of a sign for the triad gives the system the capacity to represent relations that have the character of furcations. The use of a sign for the vertice in standard graph theory packs several ideas into one conception that, for some purposes, might be better kept separate.

[Zalamea Mathematics, 2012] Fernando Zalamea, *Synthetic philosophy of contemporary mathematics*. Falmouth, UK, 2012.